\input Tex-document.sty

\def\Z{{\bf Z}}

\def\C{{\bf C}}
\def\Q{{\bf Q}}

\def\G_m{{\bf G}_m}

\def\Zbar{{\overline{\bf Z}}_{\ell}}
\def\Qbar{{\overline{\bf Q}}_{\ell}}
\def\Fbar{{\overline{\bf F}}_{\ell}}

\def\F{{\bf F}}

\let\G\Gamma

\def\Ad{\mathop{\rm Ad}\nolimits}

\def\Gal{\mathop{\rm Gal}\nolimits}

\def\Hom{\mathop{\rm Hom}\nolimits}

\def\Ad{\mathop{\rm Ad}\nolimits}

\def\id{\mathop{\rm id}\nolimits}
\def\ind{\mathop{\rm ind}\nolimits}

\def\End{\mathop{\rm End}\nolimits}

\def\Lie{\mathop{\rm Lie}\nolimits}

\def\cite{{\bf [V]}}

\def\Lie{\mathop{\rm Lie}\nolimits}
\def\Ad{\mathop{\rm Ad} \nolimits}

\def\mod{\mathop{\rm mod} \nolimits}
\def\Mod{\mathop{\rm Mod} \nolimits}
\def\id{\mathop{\rm id} \nolimits}

\def\lim{\mathop{\rm lim} \nolimits}

\bigskipamount=0.5\bigskipamount
\parskip=0.5\parskip
\hfuzz=3pt
\overfullrule=0mm

\font\boldmathHuge = cmmib10 scaled 1728


\def\firstpage{667}

\pageno=667

\def\shorttitle{Modular Representations}
\def\shortauthor{M.-F. Vign\'eras}

\title{\centerline{Modular Representations of {\boldmathHuge p}-adic Groups}
\centerline{and of Affine Hecke Algebras}}

\author{Marie-France Vign\'eras\footnote{\eightrm *}{\eightrm Institut \ de \ Math\'ematiques \ de \ Jussieu, \
Universit\'e \ de \ Paris \ 7, \ France. \ E-mail: vigneras@math.jussieu.fr}}

\vskip 7mm

\centerline{\boldnormal Abstract}

\vskip 4.5mm

{\narrower \ninepoint \smallskip
I will survey some   results in the theory of modular
representations of a reductive $p$-adic group,
 in  positive  characteristic
 $\ell \neq p$ and   $\ell=p$.

\vskip 4.5mm



}

\vskip 10mm

 {\bf Introduction\ }
   The congruences
between automorphic forms and their applications to number theory
are a motivation to study the smooth representations of a reductive
$p$-adic group $G$  over an algebraically closed field  $R$ of
{\sl any  characteristic}. The purpose of the talk is
to  give a survey of some aspects of the theory of
$R$-representations of
$G$. In
    positive  characteristic, most results are due
to the author; when proofs are available in the
litterature (some of them are not !), references will be given.

A prominent role is played by the   unipotent block which
contains the trivial representation. There is a finite list of
  types, such that the
irreducible   representations of the unipotent block are
characterized by the property that they contain a unique type of the
list.
The   types define functors
from the   $R$-representations of $G$ to the right modules  over
generalized affine Hecke algebras over $R$ with different parameters;
in positive characteristic
  $\ell$, the
  parameters  are $0$ when $\ell=p$,  and roots of unity when $\ell
\neq p$.

 In characteristic $0$ or $\ell\neq p$, for a $p$-adic linear group,
there is a Deligne-Langlands correspondence  for   irreducible
representations; the irreducible in
the   unipotent block are annihilated by  a canonical ideal
$J$;  the  category of
representations annihilated by $J$ is Morita equivalent to the
affine Schur algebra, and the unipotent  block  is  annihilated
by a finite power
$J^k$.

New phenomena appear when $\ell=p$, as
the  supersingular representations discovered by
Barthel-Livne and classified by Ch. Breuil for $ GL(2,\Q_p)$. The
modules  for the    affine Hecke algebras  of parameter $0$ and over
$R$ of characteristic $p$,    are more tractable than the
$R$-representations of the group, using that the center $Z$ of a
$\Z[q]$-affine Hecke algebra $H$ of parameter $q$ is a finitely
algebra and $H$ is a generated $Z$-module. The classification of the
simple  modules of the pro-$p$-Iwahori Hecke algebra of
$GL(2,F)$ suggests the possibility of a
  Deligne-Langlands correspondence in characteristic
$p$.

\bigskip
 {\bf Complex case\ }

Notation.    $\C$ is the field of complex numbers, $G=\underline
G(F)$ is the group of    rational points of a reductive  connected
group
$\underline G$ over a local non archimedean field $F$  with  residual
field of characteristic $p$ and of finite order
$q$, and   $\Mod_{\C}G$ is the category of complex smooth
representations of $G$. All
representations of $G$ will be smooth, the stabilizer of any vector
is open in $G$.
 An abelian category $C$ is  right (left)  Morita equivalent
to a ring
$A$ when     $C$ is equivalent to the category of right (left)
$A$-modules.

\bigskip  The
modules over complex affine Hecke algebras with parameter
$q$ are related by the Borel theorem to the complex
representations of reductive
$p$-adic groups.

\bigskip {\bf  Borel Theorem \ } {\sl The unipotent
block of
$\Mod_{\C}G$ is (left and right) Morita equivalent to  the complex
Hecke algebra of the affine Weyl group of
$G$  with parameter
$q$. }

\bigskip The proof  has three main steps,   in reverse
chronological order,   Bernstein   a)  [B]  [BK],   Borel   b) [Bo],
[C],  Iwahori-Matusmoto c) [IM], [M].

\bigskip a)  {\bf  (1.a.1) } {\sl
$\Mod_{\C}G$  is a product
of
 indecomposable abelian subcategories ``the blocks''}.

\bigskip The unipotent block   contains the trivial
representation.
The
representations in the unipotent block will be called unipotent,
although this term is   already used by Lusztig   in a
different   sense.

\bigskip {\bf  (1.a.2) }{\sl
The irreducible unipotent representations are the irreducible
subquotients of the representations parabolically induced from the
unramified characters of a minimal parabolic subgroup of $G$.}

\bigskip

b) Let $I$ be an Iwahori subgroup of $G$ (unique modulo
conjugation).

\bigskip   {\bf  (1.b.1) }{\sl The  category of
complex representations of
$G$ generated by their $I$-invariant vectors is   abelian,
equivalent by
 the functor
$$V\mapsto
V^I=\Hom_{\C G} (\C[I\backslash G], V)$$
 to the category $\Mod H_{\C}(G,I)$  of right modules of the
Iwahori Hecke algebra
$$ H_{\C}(G,I)=\End_{\C G} \C[I \backslash G]. $$}

{\bf  (1.b.2) }{\sl This abelian category  is   the unipotent
block}.

\bigskip c) {\bf  (1.c) }{\sl The Iwahori Hecke algebra
$H_{\C}(G,I)$ is the complex Hecke algebra of the affine
Weyl group of $G$   with parameter
$q$.}

\bigskip    The
algebra has a  very useful description called the
Bernstein decomposition [L1] [BK],   basic for the  geometric
description of Kazhdan-Lusztig [KL].

\bigskip From (1.b.1),  the irreducible unipotent
complex representations of
$G$ are in natural bijection with the simple modules of the complex
Hecke algebra
$H_{\C}(G,I)$.
 By the ``unipotent''   Deligne-Langlands correspondence, the
simple
$H_{\C}(G,I)$-modules ``correspond'' to the
 $G'$-conjugacy classes of pairs $(s,N)$, where $s \in
G'$ is semisimple, $N\in \Lie G'$ and $\Ad (s)N = qN$,
where
$ G'$ is the complex dual group of $G$ with Lie algebra $\Lie G'$.
This is known to be a bijection
when $G=GL(n,F)$ [Z] [R]. When
$G $ is  adjoint and unramified (quasi-split and split over a finite unramified
extension), it is also   known  to be
a bijection if one
  adds a third ingredient, a certain irreducible  geometric
representation $\rho$ of the component group of the simultaneous
centralizer of both $s$ and $N$ in $G'$; this is was done  by Chriss
[C], starting from the basic case where $G$ is split of connected
center treated by Kazhdan Lusztig  [KL] and  by Ginsburg [CG]\footnote
{  *}{ \eightrm Introduction page 18. Complex
representations of the absolute
   Weil-Deligne  group with semi-simple part
  trivial on the
{\sl inertia subgroup} (6.1)  are in natural bijection
with the $\ell$-adic representations of  the absolute Weil group
  trivial on the {\sl wild ramification subgroup}   for
any prime number $\ell \neq p$ [T] [D].  In the Deligne-Langlands
correspondence, one considers only the representations which are
Frobenius semi-simple.}. The adjoint and unramified case is
sufficient for many applications to automorphic forms; to my knowledge
the general case has not been done.

\bigskip According to R. Howe,
 the  complex  blocks  should be parametrized by types.  The
basic type, the trivial representation of an Iwahori subgroup,
is the type of the unipotent block.  An arbitrary block
should be right Morita equivalent to the Hecke algebra of the
corresponding  type. The  Hecke algebra of the  type should be a
generalized affine Hecke complex algebra  with different parameters
equal to positive powers of $p$.   This long program started in 1976
 is expected to be completed soon. The most important results are
those of Bushnell-Kutzko  for
$GL(n,F)$ [BK], of Morris for the description of the Hecke algebra of a
type [M], of Moy and Prasad for the definition
of unrefined types  [MP].

Conjecturally, the
classification of simple modules over complex generalized affine Hecke
algebras   and the theory of types
  will   give  the
classification of the complex irreducible  representations of the
reductive
$p$-adic groups.

\bigskip  We consider now the basic example, the general linear
$p$-adic group
$GL(n,F)$.
 The
   the complex irreducible
representations of $ GL(n,F)$ over $R$ are related by the
``semi-simple''  Deligne-Langlands correspondence   (proved by
Harris-Taylor  [HT1]  and
Henniart [He]), to the representations
of the Galois group
 $\Gal (\overline F/F)$  of a separable algebraic
closure $\overline F$ of $F$.

\bigskip {\bf Deligne-Langlands correspondence} {\sl

{\bf (1.d) } The blocks of
$\Mod_{\C} GL(n,F)$ are parametrized by the conjugacy classes of the
semi-simple $n$-dimensional complex representations
$\tau$ of the inertia group $I(\overline
F/F)$   which extend to  the Galois
group  $\Gal (\overline F/F)$.

{\bf (1.e) } The block parametrized by $\tau$ is  equivalent to the
unipotent block of a product of linear groups $G_{\tau}=GL(d_1,F_1)
\times
\ldots
\times GL(d_r,F_r)$ over  unramified extensions  $F_i$ of $F$
where
$\sum d_i [F_i:F]=n$.

{\bf (1.f) } The irreducible unipotent representations of $GL(n,F)$
are parametrized by the $GL(n,R)$-conjugacy classes of pairs $(s,N)$
where
$s\in GL(n,\C)$ is semi-simple,
$N\in M(n,\C)$ is nilpotent, and $sN=qNs$.}

 \bigskip {\bf  Modular case \ } Let $R$ be an
algebraically closed field  of any characteristic.  When the
characteristic of $R$ is
$0$, the  theory of  representations of $G$ is essentially like
the complex
 theory, and  the above    results   remain true
although some proofs need to be modified and this is not always in the
litterature. From now on, we will  consider   ``modular or
$\mod
\ell$''   representations, i.e. representations over $R$ of
 characteristic
$\ell >0$.

\bigskip {\bf  Banal primes \ } Although a reductive $p$-adic
group $G$ is infinite, it  behaves often as a finite group.
Given  a   property    of complex
representations of $G$  which has formally a
meaning  for
$\mod
\ell$ representations of $G$, one can usually prove that outside a
finite set of primes $\ell$, the  property remains valid. This set of
primes is called ``banal'' for the given property.

\bigskip For $\mod \ell$ representations  the Borel theorem
is false, because
the
$\mod \ell$ unipotent block of $GL(2,F)$ contains representations
without Iwahori invariant vectors when $q\equiv -1 \mod
\ell$ [V1].

\bigskip  {\bf (2) \   }
{\sl The Borel theorem is valid for
$\mod
\ell$ representations   when $\ell$  does not divide
the pro-order of any open compact subgroup of $G$.}

\bigskip  These primes   are banal for the three main steps in the
proof of the complex theorem.

\bigskip a) {\bf (2.a)  \ }
{\sl Any  prime    is banal
for the decomposition of
$\Mod_R G$ in blocks. }

 The complex proof
of (1.a.1) does not extend.  There is a new proof relying of the theory
of unrefined types [V5 III.6] when $\ell\neq p$.

\bigskip b)   (1.a.2), (1.b.1), (1.b.2) remain true
  because $\ell$  does not does
not divide the pro-order of the Iwahori subgroup $I$    and
[V2] [V4]:

\bigskip {\bf (2.b) } {\sl    Any irreducible cuspidal $\mod
\ell$-representation   of
$G$ is injective and projective in the category of $\mod
\ell$-representations of $G$ with a given central character when
$\ell$ is as in theorem 2. }

\bigskip c) Any prime $\ell$ is banal for the
Iwahori-Matsumoto step  because
the proofs of Iwahori-Matusmoto and of Morris are valid over $\Z$, and
for any commutative  ring
$A$, the Iwahori  Hecke
$A$-algebra
$$H_{A}(G,I)=\End_{A G} A[I\backslash G] \simeq
 H_{\Z}(G,I)\otimes_{\Z} A$$ is isomorphic to the Hecke $A$-algebra of
the affine Weyl group of
$G$ with parameter $q_A$ where $q_A$ is the natural
image of $q$ in $A$.

\bigskip The primes $\ell$ of theorem 2 are often called {\sl the banal
primes  of
$G$}  because   such primes are   banal for many properties. For
example, the category of $\mod
\ell$-representations of
$G$ with a given central character has finite cohomological dimension
[V4].    In the basic example $ GL(n,F)$,
$\ell$ is banal when $\ell \neq p$ and
 the multiplicative order of
$q$ modulo $\ell$ is $>n$.

\bigskip {\bf Limit primes \ } The set of   primes  banal for  (1.a.2),
(1.b.1)   is usually larger than the set of banal primes of $G$.
The primes of this set which are not banal will be
called, following Harris,  the {\sl limit primes} of $G$.  In the basic
example
$ GL(n,F)$,  the limit primes  $\ell$ satisfy
$q\equiv 1
\mod
\ell$ and
$\ell >n$ [V3].  For
number theoretic reasons,  the limit primes
 are quite important [DT] [Be] [HT2]. They
satisfy almost all the properties of the banal primes. For linear
groups, the limit primes are banal for the property that no cuspidal
representation is a subquotient of a proper parabolically induced
representation. This is may be true for
$G$ general.

\bigskip  Let $\Qbar$ be an algebraic closure of the field
$\Q_{\ell}$ of $\ell$-adic numbers, $\Zbar$  its ring of integers and
$\Fbar$ its residue field. The following statements
follow from the theory of types, or from the description of the center
of the category of $\mod \ell$ representations (the
Bernstein center).

\bigskip {\bf (3.1) \ } {\sl
The reduction  gives a surjective map from the isomorphism classes of
the irreducible cuspidal integral
$\Qbar$-representations  of $G$ to the irreducible cuspidal
$\Fbar$-representations of $G$, when $\ell$ is  a banal  or a
limit prime for $G$. }

\bigskip {\bf   Natural characteristic  \ } The interesting case
where the characteristic of $R$ is $p$ is not
yet understood.  There is a simplification: $R$-representations of
$G$ have non zero  vectors invariant by the
pro-p-radical $I_p$ of $I$. The irreducible
   are  quotients of
 $R[I_p\backslash G]$.

\bigskip Some calculations have been made for $GL(2,F)$  [BL] [Br]
[V9]. A direct classification of the
irreducible $R$-representations of $G=GL(2,\Q_p)$ [BL] [Br] and of
the pro-$p$-Iwahori Hecke $R$-algebra
$H_{R}(G,I_p)=\End_{RG}R[I_p\backslash G]$ (called a
$\mod p$ pro-$p$-Iwahori Hecke algebra)   shows:

\bigskip {\bf (4.1) } {\sl Suppose $R$ of characteristic $p$. The
pro-$p$-Iwahori functor gives a bijection between the irreducible
 $R$-representations of
$GL(2,\Q_p)$ and the simple right $H_{R}(G,I_p)$-
 modules.}

\bigskip This is the ``$\mod p$ simple
Borel theorem'' for the pro-$p$-Iwahori  group of
$GL(2,\Q_p)$. In particular $p$ is
banal for the simple version of (1.b.1) when $G=GL(2,\Q_p)$.
Irreducible $\mod p$ representations of $GL(2,F)$
  which are non subquotients of  parabolically induced
representations from a character of the diagonal torus are called
supersingular [BL]. There is a similar definition for
the $\mod p$  simple  supersingular
modules of the  pro-$p$-Iwahori Hecke algebra of $GL(2,F)$.

\bigskip {\bf (4.2) \ }  {\sl There is
  a natural bijection between the $\mod p$  simple  supersingular
modules of the  $\mod p $ pro-$p$-Iwahori Hecke algebra of $GL(2,F)$
and the
$\mod p$ irreducible dimension $2$ representations of the  absolute
Weil group of $F$. }

\bigskip This suggests the
existence of a $\mod p$ Deligne-Langlands correspondence. Some
computations are beeing made by R. Ollivier for $GL(3,F)$.

\bigskip We end this section with a new result on affine
Hecke algebras as in [L3], which is important for the theory of
representations modulo $p$.

\bigskip {\bf (4.3) \   } {\sl Let $H$ be an affine Hecke
$\Z[q]$-algebra of  parameter $q$ associated to a generalized affine
Weyl group $W$. Then the center $Z$ of $H$ is a finitely generated
$\Z[q]$-algebra  and $H$ is a finitely generated $Z$-module.}

\bigskip The key is to prove   that $H$ has a
$\Z[q]$-basis $(q^{k(w)}E_w)_{w\in W}$ where $(E_w)$ is a Bernstein
$\Z[q^{-1}]$-basis of $H[q^{-1}]$.
The assertion (4.3) was known when the parameter $q$ is invertible.

\bigskip {\bf Non natural characteristic }  $R$ an algebraically
closed field of  positive  characteristic $\ell\neq p$.   Any prime
$\ell
\neq p$ is banal for the ``simple
 Borel theorem''. The ``simple
Borel theorem'' is true $\mod \ell
\neq p$.

\bigskip {\bf  (5.1) \ } {\sl Suppose  $\ell\neq p$.  The
Iwahori-invariant functor gives a bijection between the irreducible
 $R$-representations of
$G$   with
$V^I\neq 0$ and the simple right
$ H_{R}(G,I)$-modules.}

\bigskip The existence of an
Haar measure on $G$ with values in $R$ implies that
$\Mod_RG$ is  left Morita equivalent to the convolution algebra
$H_R(G)$ of locally constant, compact distributions on $G$
with values in $R$. When
the pro-order of $I$ is invertible in $R$, the Haar measure on $G$
over $R$ normalized by $I$ is an idempotent of  $H_R(G)$, and
(5.1)  could have been already proved
by I. Schur [V3].
In general (5.1) follows from the
fact that
$R[I\backslash G]$ is ``almost projective''   [V5].

 \bigskip More generally, one expects that the Howe philosophy of types
remains  true for modular {\sl irreducible}
representations. Their  classification should reduce to the classification
of the simple modules for generalized affine Hecke $R$-algebras
of parameters equal to $0$ if $\ell=p$,  and to   roots of unity if
$\ell
\neq p$.
  This is known for linear groups if $\ell\neq p$
[V5]  or in characteristic $\ell=p$  for $GL(2,F)$ [V9].

\bigskip The unipotent block    is
described by a finite set
$S$ of modular types, the ``unipotent types'' [V7]. The set $S$
contains the class of the basic type $(I,\id)$. In the banal or limit
case, this is the only element of $S$.  A unipotent type $(P,\tau)$ is
the
$G$-conjugacy class of an irreducible
$R$-representation of a  parahoric subgroup $P$ of $G$,
  trivial on the pro-$p$-radical $P_p$,
cuspidal as a representation of $P/P_p$ (the group of
rational points of a finite reductive group over the residual field of
$F$).    The
isomorphism class of the compactly induced representation
$\ind_P^G \tau$ of $G$ determines the $G$-conjugacy class of
$\tau$, and conversely.
 We have
$\ind_I^G \id= R[I\backslash G]$.

\bigskip {\bf (5.2) \ Theorem \ } {\sl Suppose  $\ell\neq p$. There
exists a finite
 set $S$ of types, such that

- $\ind_P^G \tau$ is unipotent for any
$(P,\tau)\in S$,

- an
irreducible  unipotent
$R$-representation $V$ of $G$ is
 a quotient of $\ind_P^G
\tau$ for a unique $(P,\tau)\in S$, called the type of $V$,

- the map $V\mapsto \Hom_{RG} (\ind_P^G \tau, V)$ between the
irreducible quotients of $\ind_P^G \tau$ and the right
$ H_{R}(G,\tau)=\End_{RG}\ind_P^G \tau$ modules
  is a bijection. }

 \bigskip   The set $S$
has been explicitely described only when $G $ is a linear group
[V5].  In the example of
$GL(2,F)$  and  $q\equiv -1$ modulo  $ \ell$,
the set $S$ has two elements, the basic class and the class of
$(GL(2,O_F),\tau)$ where
$\tau$ is the cuspidal representation of dimension $q-1$ contained in
the reduction modulo
$\ell$ of the Steinberg representation of the finite group
$GL(2,\F_q)$.

The   Hecke algebra $ H_{R}(G,\tau)$ of  the type $(P,\tau)$
could  probably be
  described  a generalized affine Hecke $R$-algebra with
different parameters (complex case  [M] [L2], modular case
for a finite group  [GHM]).

\bigskip {\bf   The linear group in the non natural characteristic }
We consider the basic example
$G=GL(n,F)$ and $R$ an algebraically closed field of  positive
characteristic $\ell\neq p$.

\bigskip {\bf (6.1)  \ } {\sl Any prime $\ell \neq p$ is banal
for the Deligne-Langlands correspondence.}

\bigskip
This means that  (1.d) (1.e) (1.f) remain true when
$\C$  is replaced by $R$. The proof is done by constructing
congruences between automorphic representations for unitary groups
of compact type [V6].

\bigskip The unipotent
block  is partially described by the {\sl affine
Schur algebra}
$$S_R(G,I)=\End_{RG}V, \ \ \ V=\oplus_{P\supset I}\ind_P^G \id ,$$
 which is the ring of endomorphisms   of the direct  sum of the
representations of
$G$ compactly induced from the trivial representation  of the parahoric
subgroups
$P$  containing  the Iwahori subgroup $I$.
The functor of $I$-invariants gives an
isomorphism from the endomorphism ring of the $RG$-module $V$ to the
endomorphism ring of the right $H_{R}(G,I)$-module $V^I$ and the
$(S_R(G,I), H_{R}(G,I))$ module $V^I$  satisfies
the double centralizer property [V8].

\bigskip {\bf (6.2) \ } {\sl $\End_{H_{R}(G,I)}V^I=S_{R}(G,I),
\
\
\End_{S_{R}(G,I)}V^I=H_{R}(G,I).$}

\bigskip In the complex case, the affine Schur algebra $S_{\C}(G,I)$
is isomorphic to an algebra already defined R.M. Green [Gr]:
A complex  {\sl affine  quantum linear
group}
$\hat U(gl(n,q))$ has a remarkable representation
$W$ of countable dimension such that the
 tensor space $ W^{\otimes n}$ satisfies the double centralizer
property
 $$\End_{\hat S( n,q)}W^{\otimes n}= \hat H( n,q), \ \ \
\End_{\hat H( n,q)}W^{\otimes n}= \hat S( n,q)$$
where $\hat S( n,q)$ is the image of the action of $\hat
U(gl(n,q))$ in $ W^{\otimes n}$.
{\sl The algebras $\hat S( n,q)$ and
$\hat H( n,q)$ are respectively
 isomorphic  to $S_{\C}(G,I)$  and
$H_{\C}(G,I)$; the bimodules $ W^{\otimes n}$ and
$V^I$ are isomorphic.}

\bigskip  Let $J$ be
the   annihilator of
$R[I\backslash G]$ in the global Hecke algebra $H_R(G)$.

\bigskip  {\bf (6.3) \ Theorem } {\sl Suppose  $\ell \neq p$.

There exists an integer
$k>0$ such that the unipotent block of $\Mod_R G$ is the set of
$R$-representations of $G$  annihilated by $J^k$.

An irreducible  representation  of $G$ is unipotent if and only
if it is a   subquotient  of
$R[I\backslash G]$, if and only if it is annihilated by $J$.

The abelian subcategory of  representations of $G$
annihilated by
$J$  is Morita equivalent to the affine Schur algebra
$S_R(G,I)$.}

\bigskip     This  generalizes the Borel theorem to $\mod \ell$
representations when
$G$ is  a linear group.
The affine
Schur algebra exists  and the double centralizer
property (6.2)
 is true for a general reductive $p$-adic group $G$;  in the
banal case,  the affine Schur algebra is Morita equivalent to the
affine Hecke algebra.

\bigskip {\bf  Integral structures  \ } Let $\ell$ be any prime
number. There are two notions of integrality for an
admissible
$\Qbar$-representation
$V$ of
$G$,  $\dim V^K<\infty$ for all
open compact subgroups $K$ of $G$, which coincide when $\ell\neq p$
[V3].  One says that
$V$ is integral if $V$ contains a $G$-stable
$\Zbar$-submodule generated by a $\Qbar$-basis  of $V$, and $V$
is locally integral if the $H_{\Qbar}(G,K)$-module $V^K$ is integral,
i.e. contains
  a
$H_{\Zbar}(G,K)$-submodule
$\Zbar$-generated by a $\Qbar$-basis  of $V^K$, for all $K$.

\bigskip When  $V$  is irreducible and integral,
the action of the center
$Z$ of $G$ on $V$,  the central character,   is integral, i.e. takes
values in $\Zbar$. The situation is similar for a simple integral
$H_{\Qbar}(G,I)$-module $W$. The central character is
integral,  i.e. its restriction to the center of $ H_{\Zbar}(G,I)$
takes values in
$\Zbar$.

\bigskip {\bf (7.1)  Th\'eor\`eme \ }

\nobreak{\sl a) An irreducible cuspidal $\Qbar$-representation $V$ of
$G$ is
  integral if and only if its central character is integral.

b) A simple $H_{\Qbar}(G,I)$-module is integral if and only if its
central character is integral.

c) An irreducible representation $V$ of $G$ with $V^I\neq 0$
is locally integral if and only if $V^I$ is an
integral $H_{\Qbar}(G,I)$-module}.

\bigskip
The assertion   b)   results from (4.3). For a)  [V3]. For $\ell =p$,
    c)   is due to J.-F. Dat, using its    theory of $\ell$-adic
analysis [D].

\bigskip A general irreducible
$\Qbar$-representation $V$ of $G$ is contained in a
parabolically induced representation of an  irreducible cuspidal
 representation $W$ of a Levi subgroup of $G$.
 If $W$ is integral then $V$ is integral, but the converse is false
when $\ell =p$. When
$\ell
\neq p$, the converse is proved for classical groups by Dat
  using
 results of Moeglin (there is a gap in the
``proof'' of the converse in [V3]).

\bigskip {\bf (7.2) \  Brauer-Nesbitt principle  \  } [V3][V11] \
 {\sl  When
$\ell \neq p$, the integral structures
$L$ of an irreducible  $\Qbar$-representation  of $G$  are
$\Zbar G$-finitely generated (hence  commensurable)  and their
reduction
$L\otimes\Fbar$ are finite length
$\Fbar$-representations  of $G$ with the same semi-simplication (modulo
isomorphism).}

\bigskip  When $\ell=p$, this is false. An integral cuspidal
irreducible
$\overline\Q_p$-representation $V$ of
$G$  embeds in $\overline\Q_p[\Gamma \backslash
G]$, for any  discrete
co-compact-mod-center  subgroup $\Gamma$
 of $G$,  and has a natural integral structure with an admissible
reduction [V10]. When the  theory of types is known,   $V$
is induced from an open compact-mod-center subgroup, hence has an
integral  structure  with a non admissible reduction, which is not
commensurable with the first one.

\parindent 12mm

\head {\boldLARGE References}\def\itembr[#1]{\item{[#1]}}

\itembr[BL]Barthel L., Livne R., \
Modular representations of $GL_2$ of a local field: the ordinary
unramified case, J. of Number Theory 55, 1995, 1--27. Irreducible
modular representations of $GL_2$ of a local field, Duke Math. J.
75, 1994, 261--292.

\itembr[Be]Bella\"iche Jo\"el, \ Congruences endoscopiques et
repr\'esentations galoisiennes, Th\`ese Orsay 2002.

\itembr[B]Bernstein J.N.,  \  Le ``centre'' de Bernstein. Dans
J.N. Bernstein, P. Deligne,
 D. Kazhdan, M.-F.Vign\'eras, Repr\'esentations des groupes r\'eductifs sur un corps
 local, Travaux en cours. Hermann Paris 1984.

\itembr [Bo]Borel Armand, \ Admissible representations of a
semisimple group over a local field with vectors fixed under an
Iwahori subgroup,  Invent. Math. 35,   (1976), 233--259.

\itembr[Br]Breuil Christophe, \  Sur quelques repr\'esentations
modulaires et $p$-adiques de $GL(2,\Q_p)$ I, II, Preprints 2001.

\itembr[BK]Bushnell Colin, Kutzko Phillip, \ The admissible dual of
$GL(N)$ via compact open subgroups, Annals of Math. Studies, Princeton
Univeristy Press, 129 (1993). Smooth representations of reductive
p-adic groups: Structure theory via types, Proc. London Math. Soc. 77,
 (1988), 582--634.

\itembr[Ca]Cartier Pierre, \ Representations of $p$-adic
groups: a survey, Proc. of Symp. in pure math. AMS XXXIII, part 1,
 1979, 111--156.

\itembr[C]Chriss, Neil A., \ The classification of
representations of unramified Hecke algebras,  Math. Nachr.  191
(1998), 19--58.

\itembr[CG]Chriss N., Ginzburg V., \ Representation theory and
complex geometry, Birkhauser 1997.

\itembr[D]Dat Jean-Francois, \ Generalized tempered
representations of
$p$-adic groups, Preprint 2002.

\itembr[DT]Diamond Fred, Taylor Richard, \ Non-optimal levels of mod
$\ell$ representations, Invent. math. 115, (1994), 435--462.

\itembr[GHM] Geck, Meinolf; Hiss, Gerhard; Malle, Gunter, \ Towards a
classification of the irreducible representations in
non-describing characteristic of a finite group of Lie type.
Math. Z.  221  (1996), no. 3, 353--386.

\itembr[Gr]Green R.M., \  The affine $q$-Schur algebra,
Journal of Algebra 215 (1999) 379--411.

\itembr[IM]Iwahori N., Matsumoto H., \ On some Bruhat
decompositions and the structure of the Hecke  rings of p-adic
Chevalley groups, Publ. Math. I.H.E.S. 25 (1965),   5--48.

\itembr[HT1]Harris Michael, Taylor Richard, \ The geometry and
cohomology of some simple Shimura varieties, Annals of mathematics
studies 151 (2001).

\itembr[HT2]Harris Michael, Taylor Richard, \ Notes on $p$-adic
uniformization and congruences, 2002.

\itembr[H]Henniart Guy, \ Une preuve simple des conjectures de
Langlands pour $GL_n$ sur un corps p-adique, Invent. mat. 139
(2000), 339--350.

\itembr[KL]Kazhdan D., Lusztig G., \ Proof of the
Deligne-Langlands conjecture for Hecke algebras, Invent. Math. 87,
(1987), 153--215 .

\itembr[L1]Lusztig G., \ Some examples of square integrable
representations of $p$-adic groups, \ Trans. Amer. Math. Soc. 277,
 (1983), 623--653.

\itembr[L2]Lusztig G., \ Classification of unipotent
representations of simple $p$-adic groups, International Mathematics
Research Notices
 No11, (1995), 517--589 .

\itembr[L3]Lusztig G., \ Representations of affine Hecke algebras,
Soc. Math. de France, Ast\'erisque 171-171 (1989), 73--84.

\itembr[M]Morris Lawrence, \ Tamely ramified intertwining
algebras, Invent. math. 114,  (1993),1-54. Tamely ramified
supercuspidal representations,  Ann. Sci. École Norm. Sup. (4)  29
  no. 5, (1996), 639--667.

\itembr[MP] Moy, Allen; Prasad, Gopal. Jacquet functors and unrefined
minimal $K$-types,  Comment. Math. Helv.  71  (1996),  no. 1,
98-121. Unrefined minimal $K$-types for $p$-adic groups. Invent.
Math. 116 (1994), no. 1-3, 393--408.

\itembr[R]Rogawski John, \ On modules over the Hecke algebra
of a $p$-adic group, Invent. math. 79, (1985),  443--465.

\itembr[V1]Vign\'eras M.-F., \ Repr\'esentations modulaires
 de $GL(2,F)$ en caract\'eristique $\l$, $F$
        corps p-adique, $p\neq\l$, Compositio Mathematica 72 (1989), 33--66.
Erratum,  Compositio Mathematica 101,(1996), 109--113.

\itembr[V2]Vign\'eras M.-F., \ Banal Characteristic for
Reductive p-adic Groups,  J. of Number Theory Vol.47, Number 3,
1994, 378--397.

\itembr[V3]Vign\'eras M.-F., \ Repr\'esentations
$\l$-modulaires d'un groupes r\'eductif p-adique avec $\l \neq p$,
Birkhauser Progress in Math.137 (1996).

\itembr[V4]Vign\'eras M.-F., \ Cohomology of sheaves on the
building and\break $R$-representations,  Inventiones Mathematicae 127, 1997, 349--373.

\itembr[V5]Vign\'eras M.-F., \ Induced representations of
reductive p-adic groups in characteristic $\l \neq p$,
  Selecta Mathematica New Series 4 (1998) 549--623.

\itembr[V6]Vign\'eras M.-F., \ Correspondance locale de
Langlands semi-simple pour $GL(n,F)$ modulo $\ell\neq p$,
Inventiones  144, 2001, 197--223.

\itembr[V7]Vign\'eras M.-F., \ Irreducible modular
representations of a reductive p-adic group and simple modules for Hecke algebras,  International European\break
Congress Barcelone 2000. Birkhauser Progress in Math. 201, 117--133.

\itembr[V8]Vign\'eras M.-F., \ Schur algebra of
reductive p-adic groups I,
Institut de Math\'ematiques de Jussieu, pr\'epublication 289,  Mai
2001, To appear in
 Duke Math. Journal

\itembr[V9]Vign\'eras M.-F., \ Representations modulo p of the
p-adic group $GL(2,F)$, Institut de Math\'ematiques de Jussieu,
pr\'epublication 30,  septembre 2001.

\itembr[V10]Vign\'eras M.-F., \ Formal degree and existence of stable
arithmetic lattices of cuspidal representations of $p$-adic
reductive groups,  Invent. Math.  98 no. 3, (1989),   549--563.

\itembr[V11]Vign\'eras M.-F., \ On highest Whittaker models and
integral structures, Institut de Math\'ematiques de Jussieu,
pr\'epublication 308, septembre 2001.

\itembr[Z]Zelevinski A., \ Induced representations of
reductive p-adic groups II, Ann. scient. Ecole Norm. Sup. tome 13,
 (1980), 165--210.

\end

\font\teneurb=eurb10 \font\seveneurb=eurb7 \font\fiveeurb=eurb5
\newfam\eurbfam
\textfont\eurbfam=\teneurb \scriptfont\eurbfam=\seveneurb \scriptscriptfont\eurbfam=\fiveeurb

\font\tenmsb=msbm10 \font\sevenmsb=msbm7 \font\fivemsb=msbm5
\newfam\msbfam
\textfont\msbfam=\tenmsb \scriptfont\msbfam=\sevenmsb \scriptscriptfont\msbfam=\fivemsb

 \font\seveneuf=eufm7 \font\fiveeuf=eufm5 
\newfam\euffam
\scriptfont\euffam=\seveneuf \scriptscriptfont\euffam=\fiveeuf



\def\doublemap{\mathrel{\null _{\raise.2ex\hbox{$\textstyle\rightarrow$}}
  ^{\ellower.2ex\hbox{$\textstyle\rightarrow$}}}}

\def\date{le\space\the\day \ifcase\month\or janvier \or f\'evrier\or mars\or
avril\or mai\or juin\or juillet\or ao\^ut\or septembre\or octobre\or novembre\or d\'ecembre\fi\
{\oldstyle\the\year}}

\def\hfl#1#2{\smash{\mathop{\hbox to 12mm{\rightarrowfill}}
\limits^{\scriptstyle#1}_{\scriptstyle#2}}}

\def\"#1{\if#1i{\accent"7F\i}\else{\accent"7F#1}\fi}
\def\^#1{\if#1i{\accent"5E\i}\else{\accent"5E#1}\fi}

\def\doublemap{\mathrel{\null _{\raise.2ex\hbox{$\textstyle\rightarrow$}}
  ^{\ellower.2ex\hbox{$\textstyle\rightarrow$}}}}

\def\date{le\space\the\day \ifcase\month\or janvier \or f\'evrier\or mars\or
avril\or mai\or juin\or juillet\or ao\^ut\or septembre\or octobre\or novembre\or d\'ecembre\fi\
{\oldstyle\the\year}}